\documentclass{article}

\usepackage{delarray,verbatim,enumerate,a4wide}
\usepackage{graphicx}
\usepackage{xcolor}
\usepackage{amsmath,amsthm,amstext,amsbsy,amssymb,amsfonts,amscd}%,Mnsymbol}
\usepackage[all]{xy}
\usepackage[frenchb]{babel}
\usepackage[T1]{fontenc}
\usepackage[utf8]{inputenc}

\usepackage{relsize}

\usepackage{upgreek}
\usepackage{color}
\definecolor{vert}{rgb}{0.1,0.4,0.2}
\usepackage[colorlinks=true,linkcolor=blue,citecolor=vert]{hyperref}

\usepackage{calligra}
\DeclareFontShape{T1}{calligra}{m}{n}{<->s*[0.95]callig15}{}
\DeclareMathAlphabet{\mathscr}{T1}{calligra}{m}{n}

\newtheorem{Th}{Théorème}[]
\newtheorem{Lem}[Th]{Lemme}
\newtheorem{Prop}[Th]{Proposition}
\newtheorem{Cor}[Th]{Corollaire}

\newtheorem{Sco}[Th]{Scolie}

\newtheorem{PDef}[Th]{Proposition \& Définition}
\newtheorem{DProp}[Th]{Définition \& Proposition}

\newtheorem{Exe}[Th]{Exemple}

\def\Preuve{\smallskip\noindent {\it Preuve.~}}

\font\teneufm=eufm10
\font\seveneufm=eufm7
\font\fiveeufm=eufm5
\newfam\gothfam
\textfont\gothfam=\teneufm \scriptfont\gothfam=\seveneufm
\scriptscriptfont\gothfam=\fiveeufm

		\def\QQ{\mathbb Q}	
\def\NN{\mathbb N}	\def\ZZ{\mathbb Z}
\def\F2{\mathbb{F}_2}	\def\Z2{\mathbb{Z}_2}		
\def\Zl{\mathbb{Z}_\ell} 		

 		\def\P{\mathcal  P}			\def\F{\mathcal  F}
\def\J{\mathcal  J}  		\def\R{\mathcal  R}	\def\D{\mathcal  D}
 	  	\def\Cl{\mathcal  C \!\ell}
\def\E{\mathcal  E}		\def\T{\mathcal  T}			\def\X{\mathcal X}

	\def\p{{\mathfrak p}}			\def\a{{\mathfrak a}}
\def\L{{\mathfrak L}}		\def\l{{\mathfrak l}}

	\def\div{\operatorname{div}}
	\def\deg{\operatorname{deg}}
\def\Gal{\operatorname{Gal}}		\def\Rad{\operatorname{Rad}}
\def\Ker{\operatorname{Ker}}	\def\Coker{\operatorname{Coker}}	\def\Hom{\operatorname{Hom}}

\makeatletter
\newcommand*\wt[2][0.2ex]{%
        \begingroup
        \mathchoice{\wt@helper{#1}{#2}{\displaystyle}{\textfont}}
                   {\wt@helper{#1}{#2}{\textstyle}{\textfont}}
                   {\wt@helper{#1}{#2}{\scriptstyle}{\scriptfont}}
                   {\wt@helper{#1}{#2}{\scriptscriptstyle}{\scriptscriptfont}}%
        \endgroup
        #2%
}
\newcommand*\wt@helper[4]{%
        \def\currentfont{\the#41}%
        \def\currentskewchar{\char\the\skewchar\currentfont}%
        \setbox\tw@\hbox{\currentfont$#2$\currentskewchar}%
        \dimen@ii\wd\tw@
        \setbox\tw@\hbox{\currentfont$#2${}\currentskewchar}%
        \advance\dimen@ii-\wd\tw@
        \rlap{\raisebox{-#1}{$\m@th#3\kern\dimen@ii\widetilde{\phantom{#2}}$}}%
}
\makeatother

\def\wE{\,\wt[0.2ex]{\!\mathcal E}}		
	\def\wCl{\wt[0.1ex]{\mathcal C\!\ell}}

\newcommand{\Bmu}{\mbox{$\raisebox{-0.59ex}
 {$l$}\hspace{-0.18em}\mu\hspace{-0.88em}\raisebox{-0.98ex}{\scalebox{2}
 {$\color{white}.$}}\hspace{-0.416em}\raisebox{+0.88ex}
 {$\color{white}.$}\hspace{0.46em}$}{}}
\newcolumntype{x}[1]{>{\centering\hspace{0pt}}p{#1}}
\begin{document}

\title{\LARGE\bf  Sur la capitulation pour le module de Bertrandias-Payan}

\author{ Jean-François {\sc Jaulent} }
\date{}
\maketitle
\bigskip\bigskip\bigskip\bigskip

{\small
\noindent{\bf Résumé.}
Nous déterminons l'ordre du sous-groupe de capitulation pour le module de Bertrandias-Payan dans une $\ell$-extension arbitraire de corps de nombres qui satisfait la conjecture de Leopoldt. Nous relions en particulier la question de sa trivialisation au problème des tours localement cyclotomiques. }

\smallskip

{\small
\noindent{\bf Abstract.}
We compute the capitulation kernel  for the module of Bertrandias-Payan in an arbitrary $\ell$-extension of number fields which satisfies the Leopoldt conjecture. As a consequence we relate the existence of extensions with trivial such module to the classical problem of locally cyclotomic towers.}
\bigskip\bigskip

%%%%%%%%%%%%%%%%%%%%%%%%%%%%%%%%%%%%%%%%%%%

%%%%%%%%%%%%%%%%%%%%%%%%%%%%%%%%%%%%%%%%%%%
\bigskip

\noindent{\bf  Introduction}\medskip

%%%%%%%%%%%%%%%%%%%%%%%%%%%%%%%%%%%%%%%%%%%

Le compositum $K^{bp}$ des $\ell$-extensions cycliques d'un corps de nombres $K$ qui sont localement plongeables dans une $\Zl$-extension de ce corps a été introduit par Bertrandias et Payan dans \cite{BP} dans le but d'exhiber une condition suffisante de la conjecture de Leopoldt. Le sous-groupe de $\Zl$-torsion du groupe de Galois $\Gal(K^{bp}/K)$, appelé par Nguyen Quang Do {\em module de Bertrandias-Payan}, mesure, en effet, l'obstruction à ce qu'une pro-$\ell$-extension abélienne de $K$ soit {\em globalement} plongeable dans un produit de $\Zl$-extensions dès lors qu'elle l'est {\em localement}.\par

Il est naturel de se demander ce que devient cette obstruction lorsqu'on grossit le corps de base, plus précisément lorsqu'on remplace $K$ par une $\ell$-extension $L$.\smallskip

C'est précisément de ce problème que nous traitons ici.\smallskip

Nous adoptons dans cette étude le point de vue infinitésimal introduit dans \cite{J9}: toute $\ell$-extension $\Zl$-plongeable étant $\ell$-ramifiée (i.e. non ramifiée en dehors des places divisant $\ell$), l'extension $K^{bp}/K$ est contenue dans la pro-$\ell$-extension abélienne $\ell$-ramifiée maximale $M/K$. Et son groupe de Galois est donc un quotient de celui $\X_K$ de $M/K$. Or ce dernier s'interprète comme le quotient du $\Zl$-module libre $\,\D'_K$ construit sur les idéaux premiers de $K$ qui ne divisent pas $\ell$ par le sous-groupe des diviseurs principaux engendrés par les éléments d'image locale triviale aux places divisant $\ell$, et que nous appelons pour cela {\em infinitésimaux}.\par
De façon semblable, le groupe $\Gal(K^{bp}/K)$ s'écrit comme le quotient $\,\Cl_K^{\,bp}=\D'_K/ P_K^{\,bp}$ de $\,\D'_K$ par le sous-groupe des diviseurs principaux engendrés par les éléments dont une puissance est infinitésimale, et que nous appelons par conséquent {\em pseudo-infinitésimaux}.\smallskip

 Dans ce contexte, le noyau du morphisme d'extension $\,\Cl_K^{\,bp} \rightarrow \Cl_L^{\,bp}$ est formé par les classes des diviseurs de $K$ qui se principalisent dans $L$; et il s'analyse donc comme une {\em capitulation}. Le résultat principal de cette étude (Th. \ref{Capitulation} infra) consiste à le déterminer avec précision.\medskip

\noindent {\bf Nota}. Quoique personnel, ce travail a été mené en étroite concertation avec G. Gras et T. Nguyen Quang Do que je remercie ici. Une synthèse présentant les différents points de vue et résultats suivra prochainement. Mes remerciements vont aussi à Karim Belabas, qui a effectué avec {\sc pari} la plupart des calculs numériques sur les classes logarithmiques présentés dans cette étude.

%%%%%%%%%%%%%%%%%%%%%%%%%%%%%%%%%%%%%%%%%%%
\bigskip

\noindent{\bf  Convention}\medskip

%%%%%%%%%%%%%%%%%%%%%%%%%%%%%%%%%%%%%%%%%%%

Dans tout ce qui suit $\ell$ désigne un nombre premier impair, $K$ un corps de nombres arbitraire et $L$ une $\ell$-extension de $K$ de groupe de Galois $G=\Gal(L/K)$.

%%%%%%%%%%%%%%%%%%%%%%%%%%%%%%%%%%%%%%%%%%%
\newpage

\noindent{\bf 1. Bref rappel sur les classes infinitésimales}\medskip

%%%%%%%%%%%%%%%%%%%%%%%%%%%%%%%%%%%%%%%%%%%

Le groupe de Galois $\mathcal X_K=\Gal(M/K)$ de la pro-$\ell$-extension abélienne maximale $M$ d'un corps de nombres $K$ est interprété par la théorie du corps de classes comme groupe des {\em classes infinitésimales} de $K$ (cf. \cite{J9}), i.e. comme le quotient
$$
\Cl^{\,\infty}_K=\mathcal D'_K/\P_K^\infty
$$
du $\Zl$-module $\mathcal D'_K=\Zl\otimes_\ZZ Id'_K$ construit sur les diviseurs étrangers à $\ell$ par le sous-module des diviseurs principaux engendrés par les éléments infinitésimaux: le {\em sous-groupe infinitésimal} est le noyau $\R_K^\infty$ du morphisme de semi-localisation $s_\ell\;:\;\R_K=\Zl\otimes_\ZZ K^\times\rightarrow\R_{K_\ell}=\prod_{\l|\ell}\R_{K_\l}$, où $\R_{K_\l}=\varprojlim K^\times/K^{\times\ell^n}$ désigne le compactifié $\ell$-adique du groupe multiplicatif du complété de $K$ en la place $\l$.

Le groupe d'unités correspondant à $\Cl_K^{\,\infty}$, autrement dit le noyau de l'épimorphisme $\R_K^\infty\rightarrow\P_K^\infty$ est le groupe $\E_K^\infty$ des {\em unités infinitésimales} de $K$; sous la conjecture de Leopoldt (pour $\ell$ et pour $K$), il est trivial (cf. e.g. \cite{J31}, \S 2.2).\medskip

Sous cette même conjecture, il suit de là qu'il n'y a pas de capitulation  infinitésimale dans une $\ell$-extension $L/K$ de corps de nombres de groupe $G$, autrement dit que le morphisme naturel $j^\infty_{L/K}:\;\Cl^{\,\infty}_K \rightarrow \Cl^{\,\infty\,G}_L$ induit par l'extension des idéaux est toujours injectif. En effet, la suite exacte des classes infinitésimales ambiges s'écrit (cf. \cite{J9}, Th. 2):
$$
1 \rightarrow Cap^\infty_{L/K} \rightarrow \P^{\infty \,G}_L/\P_K^\infty \rightarrow  \D'_L{}^G/\D_K'  \rightarrow  \Cl^{\infty\, G}_L/j^\infty_{L/K}(\Cl_K^\infty)  \rightarrow  H^1(G, \P_L^\infty)      \rightarrow       1
$$
Or, sous la conjecture de Leopoldt  (pour $\ell$ et pour $L$), on a $\P^\infty_L=\R^\infty_L$ (par $\E_L^\infty=1$), donc $\P^{\infty\, G}_L=\R^{\infty\, G}_L =\R^\infty_K= \P_K^\infty$, i.e. $\P^{\infty\, G}_L/\P_K^\infty=1$ et, par suite: $\Ker j^\infty_{L/K}= Cap^\infty_{L/K}=1$. Ainsi:

\begin{Th}\label{Injectif}
 Dans une $\ell$-extension $L/K$ de corps de nombres, le morphisme d'extension entre classes infinitésimales
$\Cl^{\,\infty}_K \to \Cl^{\,\infty}_L$ est  injectif sous la conjecture de Leopoldt pour $\ell$ dans $L$. En particulier, le sous-module de torsion $\,\T^\infty_K$ de $\Cl^{\,\infty}_K$ s'identifie à un sous-module du sous-module des points fixes ${\,\T^\infty_L}^G$ de celui de $\Cl^{\,\infty}_L$.
\end{Th}

Ce résultat, aujourd'hui classique, a été prouvé initialement dans \cite{Gr2} (Th. I.1), puis redonné dans divers contextes, notamment par \cite{J9}, \cite{J18}, \cite{Ng1}. Une preuve figure également dans la monographie de Gras (\cite{Gr1}, Th. IV.2.1. Aucun de ces travaux cependant n'est mentionné dans \cite{Seo3, Seo4}.\medskip

\noindent {\bf Nota}. Par montée et descente, on voit qu'une classe $a$ qui capitule dans une extension $L/K$ de corps de nombres vérifie banalement $a^{[L:K]}=1$. La capitulation peut donc s'étudier indifféremment sur le module $\mathcal X_K$ ou sur son sous-module de $\Zl$-torsion.\medskip

Avant d'aller plus loin et d'introduire le module de Bertrandias-Payan, donnons un résultat qui nous sera utile et qui est en quelque sorte la version infinitésimale du Théorème 90 de Hilbert:

\begin{Lem}\label{H90}
Dans une $\ell$-extension de corps de nombres $L/K$, on a toujours: $ H^1(G, \R_L^\infty)=1$ 
\end{Lem}
  
\noindent {\em Preuve.}  Partons de la suite exacte courte  qui définit le sous-groupe infinitésimal:\smallskip

\centerline{$1 \rightarrow \R_L^\infty \rightarrow \R_L \rightarrow \R_{L_\ell}=\prod_{\mathfrak L|\ell}\R_{L_\mathfrak L} \rightarrow 1$.}\smallskip

\noindent  Prenant les points fixes par $G$, nous obtenons la suite exacte longue:\smallskip

\centerline{$1 \rightarrow \R_K^\infty \rightarrow \R_K \twoheadrightarrow \R_{K_\ell} \rightarrow H^1(G, \R_L^\infty) \rightarrow H^1(G,\R_L)=1$,}\smallskip

\noindent où le groupe de cohomologie à droite est trivial (c'est le Théorème 90); ce qui donne le résultat.

\begin{Sco}
Dans une $\ell$-extension $\ell$-ramifiée de corps de nombres $L/K$ qui vérifie la conjecture de Leopoldt (pour le premier $\ell$), le sous-module de torsion $\,\T^\infty_K$ de $\Cl^{\,\infty}_K$ s'identifie au sous-module des points fixes ${\,\T^\infty_L}^G$ de celui de $\Cl^{\,\infty}_L$:\smallskip

\centerline{${\T^{\infty\,G}_L} = j^\infty_{L/K}(\T_K^\infty) \simeq \T_K^\infty$.}
\end{Sco}

%%%%%%%%%%%%%%%%%%%%%%%%%%%%%%%%%%%%%%%%%%%%%%%%%%%%%%%%%%%%%%%%%%%%%%%
%%%%%%%%%%%%%%%%%%%%%%%%%%%%%%%%%%%%%%%%%%%%%%%%%%%%%%%%%%%%%%%%%%%%%%
\bigskip

\noindent{\bf 2. Module de Bertrandias-Payan et classes pseudo-infinitésimales}\medskip

%%%%%%%%%%%%%%%%%%%%%%%%%%%%%%%%%%%%%%%%%%%%%%%%%%%%%%%%%%%%%%%%%%%%%%%
%%%%%%%%%%%%%%%%%%%%%%%%%%%%%%%%%%%%%%%%%%%%%%%%%%%%%%%%%%%%%%%%%%%%%%%

Le groupe de Galois $\Gal(K^{bp}/K)$ du compositum $K^{bp}$ des $\ell$-extensions cycliques d'un corps de nombres $K$ qui sont {\em localement} $\Zl$-plongeables (ainsi noté parce qu'il  été introduit par Bertrandias et Payan dans \cite{BP}) est décrit par la théorie $\ell$-adique du corps de classes comme le quotient\smallskip

\centerline{$\Gal(K^{bp}/K)\simeq \J_K/\prod_{\p}\mu_{K_\p}^{\phantom{lc}}\R_K$}\smallskip

\noindent du $\ell$-adifié $\J_K=\prod_\p^{\rm res}\R_{K_\p}$ du groupe des idèles de $K$ par le sous-module engendré par les $\ell$-groupes locaux $\mu_{K_\p}^{\phantom{lc}}$ de racines de l'unité et le sous-groupe $\R_K$ des idèles principaux (cf. \cite{J18} ou \cite{J31}). Pour le regarder comme groupe de classes de diviseurs, il suffit d'écrire $\J_K$ comme le produit $\J'_K\R_K$, où $\J'_K$ est le produit $\prod_{\p\nmid\ell}^{\rm res}\R_{K_\p}$ (en vertu du théorème d'approximation simultanée) et d'observer que le quotient $\J'_K/\prod_{\p\nmid\ell}\mu_{K_\p}^{\phantom{lc}}$ s'identifie canoniquement au $\ell$-groupe $\D'_K$. Il vient ainsi:

\begin{PDef}
Le groupe de Galois $\Gal(K^{bp}/K)$ du compositum des $\ell$-extensions cycliques d'un corps de nombres $K$ qui sont {\em localement} $\Zl$-plongeables s'identifie au quotient\smallskip

\centerline{$\Cl^{\,bp}_K=\mathcal D'_K/\P_K^{\,bp}$}\smallskip

\noindent du $\ell$-groupe $\D'_K$ des diviseurs de $K$ étrangers à $\ell$ par le sous-groupe $\P_K^{\,bp}$ des diviseurs principaux engendrés par les éléments  pseudo-infinitésimaux, c'est à dire les éléments de $\R_K$ qui sont localement des racines de l'unité aux places au-dessus de $\ell$.\smallskip

Nous disons que $\Cl^{\,bp}_K$ est le $\ell$-groupe des classes pseudo-infinitésimales du corps $K$ et que son sous-module de $\Zl$-torsion $\,\T_K^{\,bp}$ est le module de Bertrandias-Payan attaché à $K$.
\end{PDef}

 Le sous-groupe des principaux pseudo-infinitésimaux est ainsi donné par la suite exacte courte:\smallskip

\centerline{$1 \rightarrow \R_K^\infty \rightarrow \R_K^{\,bp} \rightarrow \mu^{\phantom{lc}}_{K_\ell}=\prod_{\l|\ell}\mu^{\phantom{lc}}_{K_\l} \rightarrow 1$.}\smallskip

Et le groupe d'unités $\E_K^{\,bp}$ correspondant à $\Cl_K^{\,bp}$ est donc le sous-groupe $\mu^{loc}_K$ des éléments de $\R_K$ qui sont localement partout des racines de l'unité; sous la conjecture de Leopoldt (pour $\ell$ et pour $K$), c'est tout simplement le $\ell$-groupe $\mu^{\phantom{l}}_K$ des racines de l'unité dans $K$ (cf. \cite{J31}, \S 2.2).\medskip

Dans ce contexte, la suite exacte des classes ambiges prend donc la forme:

\begin{Prop}\label{Ambige}
Avec les notations ci-dessus, le sous-groupe de capitulation du groupe des classes pseudo-infinitésimales dans une $\ell$-extension de corps de nombres est donné par la suite exacte:\smallskip

\centerline{$1  \rightarrow Cap^{\,bp}_{L/K} \rightarrow {\P^{\,bp\,G}_L}/\P_K^{\,bp} \rightarrow  \D'_L\!{}^G/\D_K'  \rightarrow {\Cl^{\,bp\,G}_L}/j_{L/K}(\Cl_K^{\,bp})  \rightarrow  H^1(G, \P_L^{\,bp})      \rightarrow       1$.}\smallskip

\noindent Il suit en particulier: $|Cap^{\,bp}_{L/K}| \le ({\P^{\,bp\,G}_L}:\P_K^{\,bp})$, avec égalité en cas de $\ell$-ramification.
\end{Prop}

Ce résultat fournit directement une majoration simple de la capitulation pseudo-infinitésimale: partant de la suite exacte courte qui définit le sous-module principal  pseudo-infinitésimal\smallskip

\centerline{$1  \rightarrow \mu^{loc}_L  \rightarrow \R^{\,bp}_L  \rightarrow \P^{\,bp}_L \rightarrow 1$,}\smallskip

\noindent prenant la cohomologie et comparant la suite obtenue à celle pour $K$, on obtient par le serpent:\smallskip

\centerline{$1  \rightarrow {\P_L^{\,bp\,G}} /P_K^{\,bp}  \rightarrow H^1(G,\mu^{loc}_L)  \rightarrow H^1(G,\R^{\,bp}_L)  \rightarrow \cdots$}\smallskip

En fin de compte, on a donc les deux inégalités (et l'égalité à droite, sous la conjecture de Leopoldt dans $L$):

\centerline{$|Cap^{\,bp}_{L/K}| \le ( {\P_L^{\,bp\,G}}:\P_K^{\,bp}) \le |H^1(G,\mu^{loc}_L)| =  |H^1(G,\mu^{\phantom{lc}}_L)|$.}\smallskip

En particulier (voir aussi \cite{Gr5}, Th. 2.3 et \cite{Ng2} \S2):

\begin{Th}\label{CN1Cap}
Sous la conjecture de Leopoldt il ne peut y avoir de capitulation pour les classes pseudo-infinitésimales dans une $\ell$-extension $L/K$ de corps de nombres si le corps de base $K$ ne contient pas les racines $\ell$-ièmes de l'unité; ni, s'il les contient, lorsque $L$ est un étage fini de sa $\Zl$-extension cyclotomique.
\end{Th}

\noindent {\em Preuve.}  Si $K$ ne contient pas les racines $\ell$-ièmes de l'unité, le corps $L$ ne les contient pas non plus; et, s'il les contient et que $L$ est contenue dans la $\Zl$-extension cyclotomique $K^c$ de $K$, le groupe $\mu_L$ est cohomologiquement trivial. Dans les deux cas il suit: $H^1(G,\mu^{\phantom{lc}}_L)=1$; donc: $Cap^{\,bp}_{L/K}=1$.

%\begin{Sco}
%Si l'extension $L/K$ est cyclique, le sous-groupe de capitulation $Cap^{\,bp}_{L/K}$ l'est aussi.
%\end{Sco}

%\noindent {\em Preuve.} Sous la conjecture de Leopoldt, il s'identifie, en effet, à un sous-groupe du groupe quotient  ${\P^{\,bp\,G}_L}/\P_K^{\,bp}$, qui est lui-même un sous-groupe du groupe cyclique $H^1(G,\mu^{\phantom{lc}}_L)$.

%%%%%%%%%%%%%%%%%%%%%%%%%%%%%%%%%%%%%%%%%%%%%%%%%%%%%%%%%%%%%%%%%%%%%%%
%%%%%%%%%%%%%%%%%%%%%%%%%%%%%%%%%%%%%%%%%%%%%%%%%%%%%%%%%%%%%%%%%%%%%%%
\bigskip

\noindent{\bf 3. Capitulation pseudo-infinitésimale et ramification logarithmique}\medskip

%%%%%%%%%%%%%%%%%%%%%%%%%%%%%%%%%%%%%%%%%%%%%%%%%%%%%%%%%%%%%%%%%%%%%%%

Rappelons qu'une $\ell$-extension de corps de nombres $L/K$ est dite {\em logarithmiquement non ramifiée} lorsqu'elle est {\em localement cyclotomique}, i.e. lorsque pour toute place $\l$ de $K$ et chaque place $\L$ de $L$ au-dessus de $\l$, l'extension locale $L_\L/K_\l$ est contenue dans la $\Zl$-extension cyclotomique locale $K_\l^c/K_\l$ (cf. \cite{J28} ou \cite{J31}).
\begin{itemize}
\item Aux places modérées (i.e. étrangères à $\ell$), la $\Zl$-extension cyclotomique $K_\l^c/K_\l$ coïncidant avec la $\Zl$-extension non ramifiée $K^{nr}_\l/K_\l$, la  non-ramification logarithmique coïncide avec la non-ramification au sens habituel. Une $\ell$-extension logarithmiquement non ramifiée est donc, en particulier, non ramifiée (au sens habituel) en dehors de $\ell$, i.e. $\ell$-ramifiée.
\item Aux places sauvages (i.e. divisant $\ell$), en revanche, la notion de $\ell$-extension logarithmiquement non ramifiée diffère substantiellement de celle de $\ell$-extension non-ramifiée.
\end{itemize}\smallskip

Ce point rappelé, revenons à notre problème et considérons une $\ell$-extension arbitraire $L/K$. Nous pouvons supposer $\mu^{\phantom{lc}}_K\ne 1$ et $L\ne K[\mu^{\phantom{lc}}_L]$ sans quoi il ne peut y avoir de capitulation pseudo-infinitésimale en vertu du Théorème \ref{CN1Cap} ci-dessus.\smallskip

Notons alors $K'= K[\mu^{\phantom{lc}}_L]$ et observons que le morphisme d'extension de $\,\Cl_K^{\,bp}$ dans $\,\Cl_{K'}^{\,bp}$ est injectif. Considérons un diviseur $\a_K$, d'ordre $\ell$ dans $\,\Cl_K^{\,bp}$, qui se principalise dans $\,\Cl_L^{\,bp}$, ce que nous écrivons $\a_L=(\alpha_L)$ dans $\D'_L$. Il suit: $(\alpha_L^\ell)=\a_L^\ell= (\alpha_K)$, toujours dans $\D'_L$, pour un certain $\alpha_K$ dans $\R_K^{\,bp}$. Ainsi $\alpha_L^\ell/\alpha_K$ est une unité pseudo-infinitésimale de $L$, i.e. (sous la conjecture de Leopoldt) une racine de l'unité $\zeta_{K'}\in\mu^{\phantom{lc}}_{K'} = \mu^{\phantom{lc}}_L$.
Il vient alors $\alpha_L^\ell = \alpha_K\zeta_{K'}$ et $\a_K$, qui est toujours d'ordre $\ell$ dans  $\,\Cl_{K'}^{\,bp}$ se principalise dans la sous-extension kummérienne $L'=K'[\sqrt[\ell]{ \alpha_K\zeta_{K'}}]$ de $L$, laquelle est logarithmiquement non-ramifiée sur $K'$ donc sur $K$, puisque l'élément $\alpha_K\zeta_{K'}$ est localement une racine de l'unité aux places au-dessus de $\ell$ et une puissance $\ell$-ième d'idéal aux places étrangères à $\ell$. Itérant ce procédé, nous obtenons (analoguement à \cite{Gr5} \S 2):

\begin{Prop}\label{LogRamification}
Soient $K$ un corps de nombres contenant les racines $\ell$-ièmes de l'unité et $L$ une $\ell$-extension de $K$, puis $L'$ le compositum des sous-extensions de $L$ qui sont logarithmiquement non-ramifiées sur $K$ (en d'autres termes la plus grande sous-extension de $L$ qui est localement cyclotomique sur $K$). Sous la conjecture de Leopoldt, la capitulation pseudo-infinitésimale dans $L/K$ se lit dans $L'/K$:

\centerline{$Cap^{\,bp}_{L/K}=Cap^{\,bp}_{L'/K}$.}
\end{Prop}

Ce point acquis, ce n'est donc pas restreindre la généralité que de supposer l'extension $L/K$ non-ramifiée aux places au-dessus de $\ell$; auquel cas la trivialité du quotient ${\D'_L}^{\!G}/\D'_K$ donne l'égalité $Cap^{\,bp}_{L/K} = {\P^{\,bp\;G}_L}/\P_K^{\,bp}$ et la capitulation pseudo-infinitésimale s'interprète ainsi  comme le noyau du morphisme entre groupes de cohomologie: $ H^1(G,\mu^{loc}_L)  \rightarrow H^1(G,\R^{\,bp}_L) $.
\smallskip

Pour calculer le groupe de droite, partons de la suite exacte qui définit $\R_L^{\,bp}$:\smallskip

\centerline{$1 \rightarrow \R_L^\infty \rightarrow \R_L^{\,bp} \rightarrow \mu^{\phantom{lc}}_{L_\ell}=\prod_{\l|\ell}\mu^{\phantom{lc}}_{L_\l} \rightarrow 1$.}\smallskip

\noindent Prenant la cohomologie et comparant la suite obtenue avec celle écrite pour $\R_K^{\,bp}$, nous obtenons:\smallskip

\centerline{$1 \rightarrow H^1(G, \R_L^\infty)  \rightarrow H^1(G, \R_L^{\,bp})  \rightarrow H^1(G, \mu_{L_\ell}) \rightarrow H^2(G, \R_L^\infty)$.}\smallskip

\noindent Et le groupe à gauche $ H^1(G, \R_L^\infty)$ est trivial, comme expliqué dans la section 1.

\medskip

Il suit de là que l'on a  $H^1(G, \R_L^{\,bp})=1$, sous la condition {\em suffisante} $H^1(G, \mu_{L_\ell})=1$, laquelle est automatiquement vérifiée dès lors que l'extension $L/K$ considérée est {\em localement} cyclotomique, c'est à dire, suivant la terminologie introduite dans \cite{J28}, {\em logarithmiquement non ramifiée}.\par

Il vient ainsi, sous la conjecture de Leopoldt, qui assure l'égalité $\mu^{loc}_L=\mu^{\phantom{lc}}_L$:

\begin{Prop}
Soit $L/K$ une $\ell$-extension de corps de nombres logarithmiquement non ramifiée et vérifiant la conjecture de Leopoldt. La capitulation pour le module de Bertrandias-Payan est alors donnée par l'isomorphisme:

\centerline{$Cap^{\,bp}_{L/K} \simeq H^1(G,\mu_{L}^{\phantom{lc}})$.}
\end{Prop}

%%%%%%%%%%%%%%%%%%%%%%%%%%%%%%%%%%%%%%%%%%%%%%%%%%%%%%%%%%%%%%%%%%%\bigskip
\bigskip

\noindent{\bf 4. Détermination de la capitulation pseudo-infinitésimale}\medskip

%%%%%%%%%%%%%%%%%%%%%%%%%%%%%%%%%%%%%%%%%%%%%%%%%%%%%%%%%%%%%%%

Récapitulant les propositions précédentes, nous obtenons la description complète de $Cap^{\,bp}_{L/K}$:

\begin{Th}\label{Capitulation}
Soient $L/K$ une $\ell$-extension arbitraire de corps de nombres vérifiant la conjecture de Leopoldt, $L'/K$ sa plus grande sous-extension qui est localement cyclotomique et $G'=\Gal(L'/K)$ son groupe de Galois. La capitulation  dans $L/K$ pour le module de Bertrandias-Payan est alors donnée par l'isomorphisme:

\centerline{$Cap^{\,bp}_{L/K} \simeq H^1(G',\mu_{L}^{\phantom{c}})$.}
\end{Th}

\begin{Sco}
Lorsque la $\ell$-extension $L/K$ est linéairement disjointe de la $\Zl$-extension cyclotomique $K^c/K$, la capitulation pour le module de Bertandias-Payan est donnée par l'isomorphisme:\smallskip

\centerline{$Cap^{\,bp}_{L/K} \simeq \Hom(G',\mu_{K}^{\phantom{c}})$.}
\end{Sco}

\noindent {\em Preuve.} L'hypothèse de disjonction linéaire se traduit, en effet, par l'égalité: $\mu_{L}^{\phantom{lc}}=\mu_{K}^{\phantom{lc}}$.\medskip

\begin{Cor}
Soient $K$ un corps de nombres contenant les racines $\ell$-ièmes de l'unité, $K^{lc}$ sa pro-$\ell$-extension abélienne localement cyclotomique maximale et $K^c$ sa sous-extension globalement cyclotomique. Soient $L$ une $\ell$-extension cyclique de $K$ vérifiant la conjecture de Leopoldt (pour le premier $\ell$), $L'=K^{lc}\cap L$ et $K'=K^c\cap L$. L'ordre de la capitulation pseudo-infinitésimale est alors le minimum de l'ordre du $\ell$-groupe de racines $\mu^{\phantom{l}}_K$et du degré de l'extension $L'/K'$:\smallskip
 
\centerline{$|Cap^{\,bp}_{L/K}|=\min\,\{|\mu^{\phantom{l}}_K|,[L':K']\}$.}
\end{Cor}

\Preuve Observons d'abord que nous avons par construction $K'=K[\mu^{\phantom{l}}_L]$ donc $\mu^{\phantom{l}}_{K'} = \mu^{\phantom{l}}_{L'} = \mu^{\phantom{l}}_{L}$. Le groupe quotient $G'=\Gal(L'/K)$ étant cyclique, le premier groupe de cohomologie de $G'$ relatif à $\mu^{\phantom{l}}_L$ s'écrit comme quotient

\centerline{ $H^1(G',\mu^{\phantom{l}}_L)\simeq{}_N\mu^{\phantom{l}}_L/\mu^{I_{G'}}_L$}\smallskip

\noindent du noyau de la norme ${}_N\mu^{\phantom{l}}_L=\{\zeta\in\mu^{\phantom{l}}_L\,|\, N_{L'/K}(\zeta)=1\}$ par l'image de l'augmentation $\mu^{I_{G'}}_L$. Le dénominateur $\mu^{I_{G'}}_L=\{\zeta\in\mu^{\phantom{l}}_L \,|\,N_{K'/K}(\zeta)=1\}\simeq\mu^{\phantom{l}}_{K'}/\mu^{\phantom{l}}_K$ est formé des racines $[K':K]$-ièmes de l'unité; et le numérateur est sa pré-image dans $\mu^{\phantom{l}}_L$ par l'application $\zeta\mapsto N_{L'/K'}(\zeta)=\zeta^{[L':K']}$. Le premier est ainsi d'ordre $(\mu^{\phantom{l}}_L:\mu^{\phantom{l}}_K)=[K':K]$; le second, d'ordre $\min\,\{|\mu^{\phantom{l}}_L|,[L':K]\}$.

\begin{Exe}
Soient $K$ un corps de nombres contenant les racines $\ell$-ièmes de l'unité et $L$ une extension cyclique de degré $\ell$ sur $K$ vérifiant la conjecture de Leopoldt (pour le premier $\ell$). Alors:
\begin{itemize}
\item[(i)] Si $L/K$ est cyclotomique (i.e. contenue dans $K^c/K$), la capitulation $Cap^{\,bp}_{L/K}$ est triviale.
\item[(ii)] Si $L/K$ est localement cyclotomique mais non globalement, on a: $|Cap^{\,bp}_{L/K}|=\ell$
\item[(iii)] Si $L/K$ n'est pas localement cyclotomique, la capitulation $Cap^{\,bp}_{L/K}$ est triviale.
\end{itemize}
\end{Exe}

\noindent{\bf Illustrations numériques}. % (suggérées par \cite{Gr5} et basées sur les tables de Thomas (cf. \cite{Th})).\smallskip
Prenons $\ell=3$, pour $K$ un corps biquadratique $\QQ[\sqrt{d},\sqrt{-3}]$ contenant les racines cubiques de l'unité; et notons $\varepsilon$ l'unité fondamentale  du sous-corps réel de $K$.\smallskip

\begin{itemize}
\item[(i)]  Pour $K$ arbitraire et $L=K[\cos(2\pi/9)]$, l'extension $L/K$ est le premier étage de la $\ZZ_3$-tour cyclotomique; il vient donc: $Cap^{\,bp}_{L/K}= 1$, quelle que soit la valeur de $d$.
\item[(ii)] Pour $d=42,105,195, 258$ et $L=K[\sqrt[3]{\varepsilon}]$, l'extension $L/K$ est logarithmiquement non ramifiée (mais non 3-décomposée) et n'est pas globalement cyclotomique; il suit: $|Cap^{\,bp}_{L/K}|= 3$.
\item[(iii)] Pour $d=142,223,229, 235$, en revanche,  et $L=K[\sqrt[3]{\varepsilon}]$, l'extension $L/K$ est logarithmiquement ramifiée aux places au-dessus de 3; il suit: $|Cap^{\,bp}_{L/K}|= 1$.
\end{itemize}
\medskip

\noindent{\bf Remarques}. Lorsque $K$ contient les racines $\ell$-ièmes de l'unité, on sait par \cite{J23}, Cor. A.2 ou \cite{J31}, \S 3.3 que la trivialité du sous-groupe de $\Zl$-torsion $\,\T^{\,bp}_K$ du module de Bertrandias-Payan $\,\Cl^{\,bp}_K$ équivaut à celle du $\ell$-groupe des classes logarithmiques $\,\wCl_K\simeq\Gal(K^{lc}/K^c)$.\par
Dans chacun des exemples $L=K[\sqrt[3]{\varepsilon}]$ ci-dessus, le groupe de classes logarithmiques de $K$ est d'ordre 3, ce qui assure la non-trivialité de $\,\T^{\,bp}_K$. La liste fait apparaitre toutes les valeurs $d<300$ pour lesquelles cette condition est remplie et $\varepsilon$ n'est pas un cube local aux places au-dessus de 3.
%Pour plus de détails sur le radical initial des $\Zl$-extensions, on pourra comparer \cite{Gr4}, \cite{He}, \cite{J23}, \cite{J55} et \cite{MN2}, \cite{Th} avec \cite{Seo1} et \cite{Seo2}, qui ne citent aucun des précédents.

%%%%%%%%%%%%%%%%%%%%%%%%%%%%%%%%%%%%%%%%%%%%%%%%%%%%%%%%%%%%%%%%%%%\bigskip
\newpage

\noindent{\bf 5. Conoyau du morphisme d'extension dans le cas cyclique localement cyclotomique}\medskip

%%%%%%%%%%%%%%%%%%%%%%%%%%%%%%%%%%%%%%%%%%%%%%%%%%%%%%%%%%%%%%%

Intéressons-nous maintenant au conoyau $C\!ocap^{\,bp}_{L/K}= \Coker j^{\,bp}_{LK}$  du morphisme d'extension $ j^{\,bp}_{LK}:\Cl^{\,bp}_K\rightarrow\Cl_L^{\,bp\;G}$ dans une $\ell$-extension localement cyclotomique $L/K$ de groupe de Galois $G$.\smallskip

D'apès la proposition \ref{Ambige}, la suite exacte des classes pseudo-infinitésimales ambiges nous donne:\smallskip

\centerline{$1 \rightarrow H^1(G, \P^{\,bp}_L) \rightarrow C\!ocap^{\,bp}_{L/K}  \rightarrow H^2(G,\mu_L^{\phantom{lc}})  \rightarrow H^2(G, \R^{\,bp}_L) \rightarrow \cdots$}\smallskip

\noindent ce qui, lorsque $G$ est cyclique, nous conduit à la suite exacte courte:\smallskip

\centerline{$1 \rightarrow H^1(G, \P^{\,bp}_L) \rightarrow C\!ocap^{\,bp}_{L/K}  \rightarrow (\mu_K^{\phantom{lc}} \cap N_{N/K}( \R^{\,bp}_L))/N_{L/K}(\mu_L^{\phantom{lc}}) \rightarrow 1$.}\smallskip

Examinons successivement les deux groupes à droite et à gauche. Il vient:

\begin{Lem}
Dans une $\ell$-extension localement cyclotomique de corps de nombres $L/K$, on a:\smallskip

\centerline{$ H^1(G, \R_L^{\,bp})=1$ .}
\end{Lem}

\noindent {\em Preuve.} Partant de la suite exacte courte qui relie éléments infinitésimaux et éléments pseudo-infinitésimaux:\smallskip

\centerline{$1 \rightarrow \R_L^\infty \rightarrow \R_L^{\,bp} \rightarrow \mu^{\phantom{lc}}_{L_\ell}=\prod_{\L|\ell}\mu^{\phantom{lc}}_{L_\L} \rightarrow 1$.}\smallskip

\noindent prenant la cohomologie et comparant la suite obtenue à celle écrite pour le corps $K$, nous obtenons par le serpent la suite exacte:\smallskip

\centerline{$1 \rightarrow H^1(G, \R^\infty_L) \rightarrow H^1(G, \R^{\,bp}_L) \rightarrow H^1(G,\mu^{\phantom{lc}}_{L_\ell}) \rightarrow \cdots$.}\smallskip

\noindent où le groupe $H^1(G, \R^\infty_L)$ est toujours nul d'après le lemme \ref{H90}; et le groupe $H^1(G,\mu^{\phantom{lc}}_{L_\ell})$ aussi dès que l'extension $L/K$ est localement cyclotomique; d'où le résultat annoncé.

\begin{Lem}
Dans une $\ell$-extension cyclique localement cyclotomique de corps de nombres $L/K$ linéairement disjointe de la $\Zl$-extension cyclotomique $K^c/K$, on a:\smallskip

\centerline{$(\mu_K^{\phantom{lc}} \cap N_{L/K}( \R^{\,bp}_L))/N_{L/K}(\mu_L^{\phantom{lc}}) = \mu_K^{\phantom{lc}}/N_{L/K}(\mu_L^{\phantom{lc}}) =  \mu_K^{\phantom{lc}}/ \mu_K^{[L:K]}$.}
\end{Lem}

\noindent {\em Preuve.} Tout le problème est de vérifier l'inclusion: $\mu_K^{\phantom{lc}} \subset  N_{L/K}( \R^{\,bp}_L)$. Partons donc d'une racine de l'unité $\zeta\in\mu_K^{\phantom{lc}}$. L'extension $L/L$ étant localement cyclotomique, son image semi-locale $s_\ell(\zeta)$ dans $\mu_{K_\ell}^{\phantom{lc}}$ est contenue dans $N_{L/K}(\mu_{L_\ell}^{\phantom{lc}})$. Il existe donc un $x_L\in\R^{\,bp}_L$ tel qu'on ait: $s_\ell(\zeta/N_{L/K}(x_L))=1$; et le quotient $\zeta/N_{L/K}(x_L))=1$ est ainsi dans $\R_K^\infty$. Par ailleurs, $\zeta$, qui est norme locale partout par hypothèse, est également une norme globale en vertu du principe de Hasse, puisque l'extension $L/K$ est supposée cyclique; de sorte que  $\zeta/N_{L/K}(x_L))=1$ est donc en fait dans $\R_K^\infty\cap N_{L/K}(\R_L)$.
Le résultat attendu résulte donc du:

\begin{Lem}
Dans une $\ell$-extension cyclique de corps de nombres $L/K$, on a:\smallskip

\centerline{$\R_K^\infty\cap N_{L/K}(\R^{\phantom{l}}_L) = N_{L/K}(\R_L^\infty)$.}
\end{Lem}

\noindent {\em Preuve.} La suite longue de cohomologie associée dans le cas cyclique à  la suite exacte canonique\smallskip

\centerline{$1 \rightarrow \R^\infty_L \rightarrow \R_L \rightarrow \R_{L_\ell}=\prod_{\L|\ell}\R_{L_\L} \rightarrow 1$}\smallskip

\noindent nous fournit la séquence exacte:\smallskip

\centerline{$ \cdots  \rightarrow H^1(G,\R_{L_\ell}) \rightarrow  H^2(G, \R^\infty_L) \rightarrow H^2(G, \R_L) \rightarrow  \cdots$.}\smallskip

\noindent où le groupe à gauche est nul par le Théorème 90 semi-local; ce qui donne l'égalité annoncée.\medskip

Récapitulant ce qui précède, nous obtenons finalement le résultat synthétique suivant:

\begin{Th}\label{Cocap}
Soit $K$ un corps de nombres contenant les racines $\ell$-ièmes de l'unité; $L/K$ une $\ell$-extension cyclique localement cyclotomique linéairement disjointe de la $\Zl$-extension cyclotomique $K^c/K$; et $G=\Gal(L/K)$ son groupe de Galois. Noyau et conoyau du morphisme d'extension $ j^{\,bp}_{LK}:\Cl^{\,bp}_K\rightarrow\Cl_L^{\,bp\;G}$ entre classes pseudo-infinitésimales sont alors donnés sous la conjecture de Leopoldt par les isomorphismes:\smallskip

\centerline{$Cap^{\,bp}_{L/K} \simeq H^1(G,\mu_{K}^{\phantom{c}}) = {}_{_{[L/K]}}\mu_{K}^{\phantom{c}} \qquad \& \qquad C\!ocap^{\,bp}_{L/K} \simeq H^2(G,\mu_{K}^{\phantom{c}}) = \mu_{K}^{\phantom{c}}/\mu_{K}^{_{[L/K]}}$.}
\end{Th}

%%%%%%%%%%%%%%%%%%%%%%%%%%%%%%%%%%%%%%%%%%%%%%%%%%%%%%%%%%%%%%%
\newpage

\noindent{\bf 6. Application aux tours localement cyclotomiques}\medskip

%%%%%%%%%%%%%%%%%%%%%%%%%%%%%%%%%%%%%%%%%%%%%%%%%%%%%%%%%%%%%%%

Il est naturel de se demander si un corps de nombres donné $K$, dont le module de Bertrandias-Payan $\,\T^{\,bp}_K$ est non-trivial, possède ou non une $\ell$-extension $L$ avec  $\,\T^{\,bp}_L=1$.\smallskip

Lorsque $K$ ne contient pas les racines $\ell$-ièmes de l'unité,  le Théorème \ref{CN1Cap} répond à la question:

\begin{Th}
Si $K$ est un corps de nombres pour lequel on a $\mu_K^{\phantom{l}}=1$ mais $\,\T^{\,bp}_K\ne 1$, toutes ses $\ell$-extensions $L$ vérifient les deux mêmes propriétés:  $\mu_L^{\phantom{l}}=1$ mais $\,\T^{\,bp}_L\ne 1$.
\end{Th}

\noindent {\em Preuve.}  D'après le Théorème \ref{CN1Cap}, en effet, le morphisme d'extension $j_{L/K}^{\,bp}:\T^{\,bp}_K\rightarrow\T^{\,bp}_L$ est injectif.\medskip

Regardons maintenant le cas où $K$ contient les racines $\ell$-ièmes de l'unité. Il vient alors:

\begin{Lem}
S'il existe une $\ell$-extension $L/K$ avec $\,\T^{\,bp}_L=1$, alors sa sous-extension localement cyclotomique (i.e. logarithmiquement non ramifiée) maximale $L'/K$ vérifie aussi $\,\T^{\,bp}_{L'}=1$.
\end{Lem}

\noindent {\em Preuve.}  D'après la Proposition \ref{LogRamification}, en effet, le morphisme d'extension $i_{L/L'}^{\,bp}:\T^{\,bp}_{L'}\rightarrow\T^{\,bp}_L$ est injectif, puisque $L/L'$ ne contient par construction aucune sous-extension logarithmiquement non ramifiée.\medskip

Il suit de là que la question de la trivialisation du module de Bertrandias-Payan peut s'appréhender dans les tours localement cyclotomiques et qu'elle se ramène donc à l'étude de  la trivialisation du $\ell$-groupe des classes logarithmiques par $\ell$-extension logarithmiquement non ramifiée, déjà conduite dans une série d'articles : \cite{J32}, \cite{J36}, \cite{J39}.  Il vient ainsi:

\begin{Th} Soit $K$ un corps de nombres contenant les racines $\ell$-ièmes de l'unité. %Sous les conjectures de Leopoldt et de Gross-Kuz'min, 
Les assertions suivantes sont équivalentes, qui impliquent les conjectures de Leopoldt et de Gross-Kuz'min:\smallskip
\begin{itemize}
\item[(i)] Il existe une $\ell$-extension $L/K$ qui admet un module de Bertrandias-Payan trivial: $\,\T^{\,bp}_L=1$.\smallskip
\item[(ii)] Il existe une $\ell$-extension localement cyclotomique $L/K$ qui vérifie: $\,\T^{\,bp}_L=1$.\smallskip
\item[(iii)]  $K$ a une $\ell$-extension localement cyclotomique $L$ logarithmiquement principale: $\wCl^{\phantom{l}}_L=1$.\smallskip
\item[(iv)] La $\ell$-tour localement cyclotomique de $K$ est finie.
\end{itemize}
\end{Th}

\Preuve Rappelons qu'il est associé dans \cite{J32} à chaque corps de nombres $K$ une $\ell$-tour localement cyclotomique $K=K_0 \subset K_1\subset K_2 \subset \cdots \subset K_n \subset \cdots K_\infty$ obtenue en prenant pour premier étage $K_1$ le compositum des sous-extensions $K'/K$ de la pro-$\ell$-extension abélienne localement cyclotomique maximale $K^{lc}/K$, qui sont linéairement disjointes de la $\Zl$-extension cyclotomique $ K^c/K$ et pour lesquelles on a $K^{lc}=K'K^c$; puis en itérant le procédé à partir de $K_1$  (cf. \cite{J32} \S3).\par

Et il est montré (cf. \cite{J32}, Th. 4) que cette tour est finie si et seulement si le corps $K$ possède une $\ell$-extension localement cyclotomique $L$ logarithmiquement principale; ce qui établit l'équivalence de $(iii)$ et de $(iv)$. Celle de $(ii)$ et $(iii)$ résulte des théorèmes de dualité cités plus haut (\cite{J23}, Cor. A.2 et \cite{J31}, \S 3.3), qui assurent que les groupes $\,\wCl^{\phantom{l}}_L$ et $\,\T_L^{bp}$ ont même rang dès lors que le corps $L$ contient les racines $\ell$-ièmes de l'unité. Enfin l'équivalence de $(i)$ et $(ii)$ est donnée par le Lemme.\par

Par ailleurs, comme expliqué dans e.g. \cite{J17}, \cite{J18}, \cite{J31}, \cite{J54} ou \cite{J55}, en présence des racines $\ell$-ièmes de l'unité la trivialité du $\ell$-groupe des classes logarithmiques entraîne la validité des conjectures de Leopoldt et de Gross-Kuz'min (dite aussi conjecture de Gross généralisée). Chacune des conditions ci-dessus implique donc la validité des conjectures de Leopoldt et de Gross-Kuz'min dans une $\ell$-extension $L$ de $K$ et donc, par hérédité, dans $K$.\medskip

\noindent{\bf Remarques}.  Pour un corps de nombres donné, il peut arriver que la tour localement cyclotomique $K_\infty/K$ soit finie comme infinie. Des critères d'infinitude, basés sur le théorème de Golod et \v Safarevi\v c, sont donnés dans \cite{J32} et \cite{J36}. Inversement quelques conséquences de la finitude sont présentées dans \cite{J39}.\par
Sur le problème de la capitulation pour le $\ell$-groupe des classes logarithmiques, voir \cite{J54}.\medskip

\noindent{\bf Exemple}. Les calculs sur {\sc pari} montrent que le corps sextique $K=\QQ[\sqrt[3]{17},\sqrt{-3}]$ a un 3-groupe des classes logarithmiques d'ordre 3, donc exactement quatre 3-extensions cycliques de degré 3 qui sont localement cyclotomiques, parmi lesquelles la cyclotomique $K[\cos(2\pi/9)]$. Pour chacune des trois autres $L$, on trouve $\,\wCl_L=1$. La 3-tour localement cyclotomique de $K$ est ainsi de hauteur 1.

%%%%%%%%%%%%%%%%%%%%%%%%%%%%%%%%%%%%%%%%%%%%%%%%%%%%%%%%%%%%%%%
\newpage

\noindent{\bf 7. Typologie du cas cyclique élémentaire}\medskip

%%%%%%%%%%%%%%%%%%%%%%%%%%%%%%%%%%%%%%%%%%%%%%%%%%%%%%%%%%%%%%%

Supposons toujours que $K$ contient le groupe $\Bmu^{\phantom{lc}}_\ell$ des racines $\ell$-ièmes de l'unité et regardons plus spécialement le cas où $L$ est une $\ell$-extension cyclique élémentaire (i.e. de degré $\ell$) localement mais non globalement cyclotomique, de sorte que nous avons en particulier $\mu^{\phantom{lc}}_L=\mu^{\phantom{lc}}_K\ne 1$.\smallskip
 
Cette situation est gouvernée par le porisme suivant (à ne pas confondre avec \cite{Gr5} \S 4.1):

\begin{DProp}
Soit $K$ un corps de nombres contenant les racines $\ell$-ièmes de l'unité; $L/K$ une $\ell$-extension cyclique élémentaire localement cyclotomique; et $\a_K\in \D'_K$ représentant une classe non-triviale de $\Cl^{bp}_K$ qui capitule dans $\Cl^{bp}_L$. Alors, sous la conjecture de Leopoldt:
\begin{itemize}
\item ou bien le diviseur $\a_K$ représente une classe non-triviale du groupe logarithmique $\wCl_K$ qui capitule dans $\wCl_L$ et nous disons que l'extension $L/K$ est de type classe;
\item ou bien le diviseur $\a_K$ est logarithmiquement principal, auquel cas on a $L=K[\sqrt[\ell]{\varepsilon_K}]$ pour une unité logarithmique $\varepsilon_K \in \wE_K$, et nous disons que l'extension $L/K$ est de type unité.
\end{itemize}
\end{DProp}

\noindent{\it Preuve.} Écrivons $\a_K=(\alpha_L)$ dans $\D'_L$ et $\a_K^\ell=(\alpha'_K)$ dans $\D'_K$ puisque les hypothèses faites imposent à la classe du diviseur $\a_K$ d'être d'ordre exactement $\ell$ dans $\Cl^{bp}_K$. Le quotient $\alpha'_K/\alpha_L^\ell$ est alors,comme déjà vu dans la section 3, une racine (locale partout donc globale) de l'unité $\zeta_K$. Remplaçant $\alpha'_K$ par $\alpha_K/\zeta_K$, nous obtenons $\alpha_K=\alpha_L^\ell$ et  $L=K[\sqrt[\ell]{\alpha_K}]$.\par
Observons maintenant que l'égalité $\a_K=(\alpha_L)$ implique $\deg_L\a_L=\deg_L\alpha_L=0$ (cf. \cite{J28}), de sorte que $\a_K$, regardé comme diviseur logarithmique de $K$, est de degré nul, ce qui permet de considérer sa classe dans le groupe logarithmique $\wCl_K$.
  Cela étant:
\begin{itemize}
\item Si elle est triviale, écrivons $\a_K=\wt{\div}(\beta_K)$ pour un $\beta_K$ de $\R_K$. Il vient alors $\alpha_K=\varepsilon_K\beta_K^\ell$ pour une unité logarithmique $\varepsilon_K$ et donc $L=K[\sqrt[\ell]{\varepsilon_K}]$, comme annoncé.
\item Inversement, si l'on a $L=K[\sqrt[\ell]{\varepsilon_K}]$ pour une unité logarithmique $\varepsilon_K$, il suit  $\alpha_K=\varepsilon_K\beta_K^\ell$ pour un $\beta_K$ de $\R_K$; puis $\a_K^\ell=\wt\div(\alpha_K)=\wt\div(\beta^\ell)$. Et $\a_k=\wt\div(\beta_K)$ est logarithmiquement principal, ce qui achève la démonstration du porisme.
\end{itemize}\medskip

\noindent{\bf Exemples}. Pour $K=\QQ[\sqrt{d},\sqrt{-3}]$ et $d=29, 62, 74, 77,\dots$, le corps biquadratique $K$ possède une unique place au-dessus de 3 et un 3-groupe des classes logarithmiques d'ordre 3. Son groupe des unités logarithmiques $\wE_K$ coïncide avec le 3-adifié $\E'_K= \ZZ_3\otimes_\ZZ E'_K$ du groupe des 3-unités et l'unité fondamentale $\varepsilon$ du sous-corps quadratique réel $\QQ[\sqrt{d}]$ est donc aussi une unité logarithmique. Le calcul montre que c'est un cube local (au dessus de 3) dans $K$. Ainsi l'extension 3-ramifiée $L=K[\sqrt[3]{\varepsilon}]$ est 3-décomposée et donc localement cyclotomique. Elle est, bien sûr, de type unité.\medskip

Pour illustrer le type classe, considérons le cas des extensions à conjugaison complexe:

\begin{Prop}
Soient $K$, contenant le groupe $\Bmu^{\phantom{lc}}_\ell$ des racines $\ell$-ièmes de l'unité, un corps à conjugaison complexe, extension quadratique totalement imaginaire d'un sous-corps $K_+$ totalement réel, et $L$ une $\ell$-extension localement cyclotomique de $K$ provenant par composition avec $K$ d'une $\ell$-extension cyclique élémentaire $L_+$ de $K_+$ autre que la cyclotomique. Dans ce cas, sous la conjecture de Gross-Kuz'min, l'extension $L/K$ est toujours de type classe.
\end{Prop}

\noindent{\it Preuve.} Notons $\tau$ la conjugaison complexe. Le nombre premier $\ell$ étant impair, chaque $\Zl$-module galoisien se décompose naturellement en ses composantes réelles et imaginaires par action des deux idempotents $e_\pm=\frac{1}{2}(1\pm\tau)$. Par hypothèse, le groupe de Galois $\Gal(L/K)\simeq\Gal(L_+/K_+)$ est réel; le radical correspondant $\Rad(L/K)\simeq \Hom(\Gal(L/K),\Bmu^{\phantom{lc}}_\ell)$ est donc imaginaire. Or, sous la conjecture de Gross-Kuz'min, la composante imaginaire du groupe des unités logarithmiques se réduit aux racines de l'unité (cf. \cite{J28}). La seule $\ell$-extension élémentaire localement cyclotomique de type unité est donc la cyclotomique, qui a été exclue.\smallskip

\noindent{\it Nota.} Dans la situation ci-dessus, la composante réelle du $\ell$-groupe des classes logarithmiques $\,\wCl_K^{e_+}\simeq\wCl_{K_+}$ est non-triviale par hypothèse. Sous la conjecture de Gross-Kuz'min, les inégalités du miroir assurent alors qu'il en est de même de la composante imaginaire $\,\wCl_K^{e_-}$, laquelle a même rang que le module de Bertrandias-Payan $\,\T^{\,bp}_{K_+}\simeq (\T^{\,bp}_K)^{\,e_+}_{\phantom{l}}$. Or, celui-ci s'injecte dans $\,\T^{\,bp}_{L_+}\simeq (\T^{\,bp}_L)^{\,e_+}_{\phantom{l}}$, puisque la capitulation $Cap^{\,bp}_{L/K} \simeq H^1(G,\mu_{L}^{\phantom{lc}})$ est imaginaire. En particulier, le $\ell$-groupe $\,\T^{\,bp}_{L_+}$ est donc toujours non-trivial dans ce cas.

%%%%%%%%%%%%%%%%%%%%%%%%%%%%%%%%%%%%%%%%%%%%%%%%%%%%%%%%%%%%%%%
\newpage

\noindent{\bf 8. Capitulation ultime dans les tours localement cyclotomiques}\medskip

%%%%%%%%%%%%%%%%%%%%%%%%%%%%%%%%%%%%%%%%%%%%%%%%%%%%%%%%%%%%%%%

Les résultats de la section 6 montrent qu'un corps de nombres $K$ dont le module de Bertrandias-Payan $\,\T^{\,bp}_K$ est non-trivial ne peut admettre de $\ell$-extension $N$ avec $\,\T^{\,bp}_N=1$ que s'il contient le groupe $\Bmu^{\phantom{lc}}_\ell$ des racines $\ell$-ièmes de l'unité et si sa $\ell$-tour localement cyclotomique $N/K$ est finie. Cette tour peut alors être décrite (non canoniquement) par empilement d'extensions cycliques élémentaires $K=K_0\subset K_1 \subset \dots \subset K_n=N$ de degrés relatifs $[K_{i+1}:K_i]=\ell$.
 L'objet de cette dernière section est de regarder plus attentivement le dernier étage $K_n/K_{n-1}$ d'une telle tour.\medskip
 
 Soit donc $K$ un corps de nombres avec $\mu^{\phantom{lc}}_K \supset \Bmu^{\phantom{lc}}_\ell$ et $\,\T^{\,bp}_K\ne 1$, possédant une $\ell$-extension $L$ de degré $\ell$ avec $\,\T^{\,bp}_L=1$. Notons $G=\Gal(L/K)$ et observons que $\,\T^{\,bp}_K$ capitule dans $L$, ce qui impose à l'extension $L/K$ d'être localement mais non globalement cyclotomique, donc de vérifier $\mu^{\phantom{lc}}_L=\mu^{\phantom{lc}}_K$, et à $\,\T_K=Cap^{bp}_{L/K}\simeq\Hom(G,\mu^{\phantom{l}}_L)$ d'être d'ordre $\ell$ (cf. \S3).\smallskip
 
 Par ailleurs, $L$ contenant les racines $\ell$-ièmes de l'unité, la trivialité de $\,\T^{\,bp}_L$ équivaut à celle du $\ell$-groupe des classes logarithmiques $\;\wCl_L$ (et implique au passage la validité dans $L$ des conjectures de Leopoldt et de Gross-Kuz'min). En d'autres termes la pro-$\ell$-extension abélienne localement cyclotomique maximale $L^{lc}$ de $L$ coïncide avec sa $\Zl$-extension cyclotomique $L^c$. Il suit $K^{lc}=L^{lc}=L^c$ et le $\ell$-groupe des classes logarithmiques $\,\wCl_K$ de $K$ est donc exactement d'ordre $\ell$.\smallskip
 
 Ce point acquis, nous pouvons observer que toutes ces propriétés se propagent à chaque étage fini $L_n/K_n$ de la $\ell$-tour cyclotomique $L^c/K^c$ construite sur l'extension $L/K$. En effet, la trivialité du $\ell$-groupe des classes logarithmiques $\,\wCl_L\simeq\Gal(L^{lc}/L^c)$ montre que la pro-$\ell$-extension (galoisienne) localement cyclotomique maximale de $L$ se réduit encore à $L^c$, de sorte que tous les  $\,\wCl_{L_n}$ sont identiquement nuls. Par dualité, il en est de même des modules $\,\T^{\,bp}_{L_n}$ de sorte que chacune des extensions $L_n/K_n$ vérifie nos hypothèses de départ. En particulier, il vient ainsi:\smallskip
 
\centerline{$|\T^{\,bp}_{K_n}|=|\wCl^{\phantom{l}}_{K_n}|=\ell \qquad \& \qquad |\T^{\,bp}_{L_n}|=|\wCl^{\phantom{l}}_{L_n}|=1$}\smallskip

\noindent à chaque étage fini $n\in\NN$ de la $\ell$-tour cyclotomique.\smallskip

Regardons maintenant le sommet $L^c/K^c$ de la tour: L'extension $L^c/K^c$ est localement cyclotomique, donc, dans le cas présent, complètement décomposée en chacune de ses places. Comme $L^c$ ne possède aucune $\ell$-extension non-triviale complètement décomposée partout, c'est précisément la $\ell$-extension (en l'occurence abélienne) complètement décomposée maximale de $K^c$. Et le groupe de Galois $\Gal(L^c/K^c)\simeq G$ est donc isomorphe à la limite projective (pour les applications normes) des $\ell$-groupes de $\ell$-classes $\,\Cl'_{K_n}$ des corps $K_n$. Autrement dit, il suit:\smallskip
 
\centerline{$\varprojlim\, \wCl^{\phantom{l}}_{K_n}= \varprojlim \,\Cl'_{K_n} \simeq G \qquad \text{ i.e.} \qquad \wCl^{\phantom{l}}_{L_n}=  \Cl'_{L_n} =1 \quad \text{pour tout } n\gg 0$;}\smallskip

\noindent puis, par un argument classique de Théorie d'Iwasawa:\smallskip
 
\centerline{$\varinjlim\, \wCl^{\phantom{l}}_{K_n}= \varinjlim \,\Cl'_{K_n} = \,\Cl'_{K^c} =1$,}\smallskip

\noindent de sorte que, contrairement aux modules $\T^{\,bp}_{K_n}$, les groupes de classes $\wCl_{K_n}$ et $\,\Cl'_{K_n}$ capitulent ultimement dans la $\ell$-tour cyclotomique à chaque étage assez grand $K_{n+1}/K_n$, avec apparition à chaque fois au niveau supérieur de classes (ambiges) qui ne proviennent pas de l'étage inférieur.\smallskip

À l'opposé, puisqu'il n'y a pas de capitulation pour les modules de Bertrandias-Payan dans la $\ell$-tour cyclotomique, tous les $\,\T^{\,bp}_{K_n}$ sont engendrés par la classe d'un même diviseur $\a_K\in\D'_K$, lequel est logarithmiquement principal pour $n\gg 0$. Il suit de là que les extensions $L_n/K_n$ sont ultimement de type unité.
En résumé, il vient:

\begin{Th}
Soit $K$ un corps de nombres contenant les racines $\ell$-ièmes de l'unité avec $\,\T^{\,bp}_K\ne1$, mais qui possède  une $\ell$-extension $L$ de degré $[L:K]=\ell$ avec $\,\T^{\,bp}_L=1$. On a alors:\smallskip
\begin{itemize}
\item $|\T^{\,bp}_{K_n}|=|\wCl^{\phantom{l}}_{K_n}|=\ell$ et $\,\T^{\,bp}_{L_n}=\wCl^{\phantom{l}}_{L_n}=1$ à chaque étage fini $L_n/K_n$ de la $\ell$-tour cyclotomique.\smallskip

\item $\,\wCl^{\phantom{l}}_{K_n}=\Cl'_{K_n}$ pour tout $n\gg 0$ avec capitulation dans $\,\wCl^{\phantom{l}}_{K_{n+1}}=\Cl'_{K_{n+1}}$; et donc $\,\Cl'_{K^c}=1$.\smallskip
\end{itemize}

 En particulier les $\ell$-extensions localement cyclotomiques $L_n/K_n$ sont ultimement de type unité.

\end{Th}

\noindent{\bf Remarque}. Le corps sextique  $K=\QQ[\sqrt[3]{17},\sqrt{-3}]$ présenté à la fin de la section 3 fournit un exemple simple de la situation décrite ci-dessus. Bien entendu, quoique totalement imaginaire, il n'est pas à conjugaison complexe, cette conjugaison étant incompatible avec le type unité (cf. \S7).
\newpage
%%%%%%%%%%%%%%%%%%%%%%%%%%%%%%%%%%%%%%%%%%%%%%%%%%%%%%%%%%%%%%%

\def\refname{\normalsize{\sc  Références}}

{\smaller

}

\bigskip\noindent
{\small
\begin{tabular}{l}
{Jean-François {\sc Jaulent}}\\
Institut de Mathématiques de Bordeaux \\
Université de {\sc Bordeaux} \& CNRS \\
351, cours de la libération\\
F-33405 {\sc Talence} Cedex\\
courriel : Jean-Francois.Jaulent@math.u-bordeaux1.fr 
\end{tabular}
}

 \end{document}